\documentclass[12pt]{article}%
\usepackage{amsmath}
\usepackage{graphicx}
\usepackage{subfigure}
\usepackage{amsfonts}
\usepackage{amssymb}%
\providecommand{\U}[1]{\protect\rule{.1in}{.1in}}
\newtheorem{theorem}{Theorem}

\newtheorem{algorithm}[theorem]{Algorithm}

\newtheorem{corollary}[theorem]{Corollary}

\newtheorem{definition}[theorem]{Definition}

\newtheorem{lemma}[theorem]{Lemma}

\newtheorem{remark}[theorem]{Remark}

\newenvironment{proof}[1][Proof]{\noindent\textbf{#1.} }{\ \rule{0.5em}{0.5em}}
\numberwithin{equation}{section}
\numberwithin{theorem}{section}
\numberwithin{table}{section}
\newcommand{\be}{\begin{equation}}
\newcommand{\ee}{\end{equation}}
\begin{document}

\title{\textbf{On The Behavior of Subgradient Projections Methods for Convex
Feasibility Problems in Euclidean Spaces}}
\author{Dan Butnariu$^{1}$, Yair Censor$^{1}$, Pini Gurfil$^{2}$ and Ethan
Hadar$^{3}\bigskip$\\$^{1}$Department of Mathematics, University of Haifa\\Mt. Carmel, Haifa 31905, Israel\\(\{dbutnaru, yair\}@math.haifa.ac.il)\bigskip\\$^{2}$Faculty of Aerospace Engineering \\Technion -- Israel Institute of Technology \\Technion City, Haifa 32000, Israel\\(pgurfil@aerodyne.technion.ac.il)\bigskip\\$^{3}$CA Labs, CA Inc.\\Yokneam 20692, Israel\\(ethan.hadar@ca.com)}
\date{April 22, 2007. Revised: December 31, 2007. Revised: February 5, 2008.}
\maketitle

\begin{abstract}
We study some methods of subgradient projections for solving a convex
feasibility problem with general (not necessarily hyperplanes or half-spaces)
convex sets in the inconsistent case and propose a strategy that controls the
relaxation parameters in a specific self-adapting manner. This strategy leaves
enough user-flexibility but gives a mathematical guarantee for the algorithm's
behavior in the inconsistent case. We present numerical results of
computational experiments that illustrate the computational advantage of the
new method.

\end{abstract}

\section{Introduction\label{sect:intro}}

In this paper we consider, in an Euclidean space framework, the method of
simultaneous subgradient projections for solving a convex feasibility problem
with general (not necessarily linear) convex sets in the consistent and
inconsistent cases. To cope with this situation, we propose two algorithmic
developments. One uses \textit{steering parameters} instead of
\textit{relaxation parameters} in the simultaneous subgradient projection
method, and the other is a strategy that controls the relaxation parameters in
a specific self-adapting manner that leaves enough user-flexibility while
yielding some mathematical guarantees for the algorithm's behavior in the
inconsistent case. For the algorithm that uses steering parameters there is
currently no mathematical theory. We present numerical results of
computational experiments that show the computational advantage of the
mathematically-founded algorithm implementing our specific relaxation
strategy. In the remainder of this section we elaborate upon the meaning of
the above-made statements.

Given $m$ closed convex subsets $Q_{1},Q_{2},\cdots,Q_{m}\subseteq R^{n}$ of
the $n$-dimensional Euclidean space, expressed as%
\begin{equation}
Q_{i}=\left\{  x\in R^{n}\mid f_{i}(x)\leq0\right\}  , \label{eq:set}%
\end{equation}
where $f_{i}:R^{n}\rightarrow R$ is a convex function, the \textit{convex
feasibility problem} (CFP) is%
\begin{equation}
\text{find a point }x^{\ast}\in Q:=\cap_{i=1}^{m}Q_{i}. \label{(CFP problem)}%
\end{equation}
As is well-known, if the sets are given in any other form then they can be
represented in the form (\ref{eq:set}) by choosing for $f_{i}$ the squared
Euclidean distance to the set. Thus, it is required to solve the system of
convex inequalities%
\begin{equation}
\text{ }f_{i}(x)\leq0,\text{ \ \ }i=1,2,\ldots,m. \label{CFP inequalities}%
\end{equation}

A fundamental question is how to approach the CFP in the inconsistent case,
when $Q=\cap_{i=1}^{m}Q_{i}=\emptyset.$ Logically, algorithms designed to
solve the CFP by finding a point $x^{\ast}\in Q$ are bound to fail and should,
therefore, not be employed. But this is not always the case. Projection
methods that are commonly used for the CFP, particularly in some very large
real-world applications (see details below) are applied to CFPs without prior
knowledge whether or not the problem is consistent. In such circumstances it
is imperative to know how would a method, that is originally known to converge
for a consistent CFP, behave if consistency is not guaranteed.

We address this question for a particular type of projection methods. In
general, sequential projection methods exhibit \textit{cyclic convergence} in
the inconsistent case. This means that the whole sequence of iterates does not
converge, but it breaks up into $m$ convergent subsequences (see Gubin, Polyak
and Raik \cite[Theorem 2]{gpr67} and Bauschke, Borwein and Lewis
\cite{bbl97}). In contrast, simultaneous projection methods generally
converge, even in the inconsistent case, to a minimizer of a proximity
function that \textquotedblleft measures\textquotedblright\ the weighted sum
of squared distances to all sets of the CFP, provided such a minimizer exists
(see Iusem and De Pierro \cite{IDP} for a local convergence proof and
Combettes \cite{Combettes} for a global one).

Therefore, there is an advantage in using simultaneous projection methods from
the point of view of convergence. Additional advantages are that (i) they are
inherently parallel already at the mathematical formulation level due to the
simultaneous nature, and (ii) they allow the user to assign weights (of
importance) to the sets of the CFP. However, a severe limitation, common to
sequential as well as simultaneous projection methods, is the need to solve an
inner-loop distance-minimization step for the calculation of the orthogonal
projection onto each individual set of the CFP. This need is alleviated only
for convex sets that are simple to project onto, such as hyperplanes or half-spaces.

A useful path to circumvent this limitation is to use subgradient projections
that rely on calculation of subgradients at the current (available) iteration
points, see Censor and Lent \cite{lc82} or \cite[Section 5.3]{CZ97}. Iusem and
Moledo \cite{moledo} studied the simultaneous projection method with
subgradient projections but only for consistent CFPs. To the best of our
knowledge, there does not exist a study of the simultaneous projection method
with subgradient projections for the inconsistent case. Our present results
are a contribution towards this goal.

The CFP is a fundamental problem in many areas of mathematics and the physical
sciences, see, e.g., Combettes \cite{c93,c96} and references therein. It has
been used to model significant real-world problems in image reconstruction
from projections, see, e.g., Herman \cite{H80}, in radiation therapy treatment
planning, see Censor, Altschuler and Powlis \cite{cap88} and Censor
\cite{cen99}, and in crystallography, see Marks, Sinkler and Landree
\cite{msl99}, to name but a few, and has been used under additional names such
as \textit{set-theoretic estimation} or the \textit{feasible set approach}. A
common approach to such problems is to use projection algorithms, see, e.g.,
Bauschke and Borwein \cite{bb96}, which employ \textit{orthogonal projections}
(i.e., nearest point mappings) onto the individual sets $Q_{i}.$ The
orthogonal projection $P_{\Omega}(z)$ of a point $z\in R^{n}$ onto a closed
convex set $\Omega\subseteq R^{n}$ is defined by%
\begin{equation}
P_{\Omega}(z):=\operatorname{argmin}\{\parallel z-x\text{ }\Vert\mid\text{
}x\in\Omega\},
\end{equation}
where, throughout this paper, $\parallel\cdot\parallel$ and $\left\langle
\cdot,\cdot\right\rangle $ denote the Euclidean norm and inner product,
respectively, in $R^{n}.$ Frequently a \textit{relaxation parameter} is
introduced so that%
\begin{equation}
P_{\Omega,\lambda}(z):=(1-\lambda)z+\lambda P_{\Omega}(z) \label{eq2}%
\end{equation}
is the \textit{relaxed projection} of $z$ onto $\Omega$ with relaxation
$\lambda.$ Many iterative projection algorithms for the CFP were developed,
see Subsection \ref{sub:proj-methods} below.

\subsection{Projection methods: Advantages and earlier
work\label{sub:proj-methods}}

The reason why the CFP is looked at from the viewpoint of projection methods
can be appreciated by the following brief comments, that we made in earlier
publications, regarding projection methods in general. Projections onto sets
are used in a variety of methods in optimization theory but not every method
that uses projections really belongs to the class of projection methods.
\textit{Projection methods} are iterative algorithms which use projections
onto sets. They rely on the general principle that projections onto the given
individual sets are easier to perform then projections onto other sets derived
from the given individual sets (intersections, image sets under some
transformation, etc.)

A projection algorithm reaches its goal, related to the whole family of sets,
by performing projections onto the individual sets. Projection
algorithms\textit{ }employ projections onto convex sets in various ways. They
may use different kinds of projections and, sometimes, even use different
types of projections within the same algorithm. They serve to solve a variety
of problems which are either of the feasibility or the optimization types.
They have different algorithmic structures, of which some are particularly
suitable for parallel computing, and they demonstrate nice convergence
properties and/or good initial behavior patterns.

Apart from theoretical interest, the main advantage of projection methods,
which makes them successful in real-world applications, is computational. They
commonly have the ability to handle huge-size problems that are beyond the
ability of more sophisticated, currently available, methods. This is so
because the building blocks of a projection algorithm are the projections onto
the given individual sets (which are easy to perform) and the algorithmic
structure is either sequential or simultaneous (or in-between).

The field of projection methods is vast and we mention here only a few recent
works that can give the reader some good starting points. Such a list
includes, among many others, the works of Crombez \cite{crombez03, crombez06},
the connection with variational inequalities, see, e.g., Aslam-Noor
\cite{noor04}, Yamada \cite{yamada} which is motivated by real-world problems
of signal processing, and the many contributions of Bauschke and Combettes,
see, e.g., Bauschke, Combettes and Kruk \cite{bauschkeetal05} and references
therein. Bauschke and Borwein \cite{bb96} and Censor and Zenios \cite[Chapter
5]{CZ97} provide reviews of the field.

Systems of linear equations, linear inequalities, or convex inequalities are
all encompassed by the CFP which has broad applicability in many areas of
mathematics and the physical and engineering sciences. These include, among
others, optimization theory (see, e.g., Eremin \cite{Eremin69}, Censor and
Lent \cite{lc82} and Chinneck \cite{chinneck03}), approximation theory (see,
e.g., Deutsch \cite{Deutsch} and references therein), image reconstruction
from projections in computerized tomography (see, e.g., Herman \cite{H80,
HM93}) and control theory (see, e.g., Boyd et al. \cite{control94}.)

Combettes \cite{comb97} and Kiwiel \cite{kiwiel97} have studied the
subgradient projection method for consistent CFPs. Their work presents more
general algorithmic steps and is formulated in Hilbert space. Some work has
already been done on detecting infeasibility with certain subgradient
projection methods by Kiwiel \cite{kiwiel96a, kiwiel96b}. However, our
approach differs from the latter in that it aims at a subgradient projection
method that \textquotedblleft will work\textquotedblright\ regardless of the
feasibility of the underlying CFP and which does not require the user to study
in advance whether or not the CFP is consistent. Further questions arise such
as that of combining our work, or the above quoted results, with Pierra's
\cite{Pierra84} product space formalism, as extended to handle inconsistent
situations by Combettes \cite{Combettes}. These questions are currently under investigation.

\section{Simultaneous subgradient projections with steering
parameters\label{sect:steering}}

Subgradient projections have been incorporated in iterative algorithms for the
solution of CFPs. The {cyclic subgradient projections} (CSP) method for the
CFP was given by Censor and Lent \cite{lc82} as follows.

\begin{algorithm}
{ \label{algo:CSP} The method of cyclic subgradient projections }(CSP).

\textbf{Initialization}: $x^{0}\in R^{n}$ is arbitrary.

\textbf{Iterative step}: Given $x^{k},$ calculate the next iterate $x^{k+1}$
by
\begin{equation}
x^{k+1}=\left\{
\begin{array}
[c]{ll}%
x^{k}-\alpha_{k}\frac{\displaystyle f_{i(k)}(x^{k})}{\displaystyle \parallel
\;t^{k}\parallel^{2}}t^{k}, & \mathrm{if}\;\;f_{i(k)}(x^{k})>0,\\
x^{k}, & \mathrm{if}\;\;f_{i(k)}(x^{k})\leq0,
\end{array}
\right.  \label{eqn:CSP}%
\end{equation}
where $t^{k}\in\partial f_{i(k)}(x^{k})$ is a subgradient of $f_{i(k)}$ at the
point $x^{k}$, and the relaxation parameters $\{\alpha_{k}\}_{k=0}^{\infty}$
are confined to an interval $\epsilon_{1}\leq\alpha_{k}\leq2-\epsilon_{2}$,
for all $k\geq0$, with some, arbitrarily small, $\epsilon_{1},\epsilon_{2}>0.$

\textbf{Control}: Denoting $I:=\{1,2,\ldots,m\},$ the sequence $\{i(k)\}_{k=0}%
^{\infty}$ is an almost cyclic control sequence on $I$. This means (see, e.g.,
\cite[Definition 5.1.1]{CZ97}) that $i(k)\in I$ for all $k\geq0$ and there
exists an integer $C\geq m$ such that, for all $k\geq0$, $I\subseteq
\{i(k+1),i(k+2),\ldots,i(k+C)\}.$
\end{algorithm}

Observe that if $t^{k}=0$, then $f_{i(k)}$ takes its minimal value at $x^{k}$,
implying, by the nonemptiness of $Q$, that $f_{i(k)}(x^{k})\leq0$, so that
$x^{k+1}=x^{k}$. Relations of the CSP method to other iterative methods for
solving the convex feasibility problem and to the relaxation method for
solving linear inequalities can be found, e.g., in \cite[Chapter 5]{CZ97}, see
also, Bauschke and Borwein \cite[Section 7]{bb96}. Since sequential projection
methods for CFPs commonly have fully-simultaneous counterparts, the
simultaneous {subgradient projections (SSP) method of Dos Santos
\cite{dossantos} and Iusem and Moledo \cite{moledo} is a natural algorithmic
development.}

\begin{algorithm}
\label{alg:ssp}{The method of simultaneous subgradient projections }(SSP){. }

\textbf{Initialization}: $x^{0}\in R^{n}$ is arbitrary.

\textbf{Iterative step}: \textbf{(i)} Given $x^{k},$ calculate, for all $i\in
I=\{1,2,\ldots,m\},$ intermediate iterates $y^{k+1,i}$ by
\begin{equation}
y^{k+1,i}=\left\{
\begin{array}
[c]{ll}%
x^{k}-\alpha_{k}\frac{\displaystyle f_{i}(x^{k})}{\displaystyle\parallel
\;t^{k}\parallel^{2}}t^{k}, & \mathrm{if}\;\;f_{i}(x^{k})>0,\\
x^{k}, & \mathrm{if}\;\;f_{i}(x^{k})\leq0,
\end{array}
\right.
\end{equation}
where $t^{k}\in\partial f_{i}(x^{k})$ is a subgradient of $f_{i}$ at the point
$x^{k}$, and the relaxation parameters $\{\alpha_{k}\}_{k=0}^{\infty}$ are
confined to an interval $\epsilon_{1}\leq\alpha_{k}\leq2-\epsilon_{2}$, for
all $k\geq0$, with some, arbitrarily small, $\epsilon_{1},\epsilon_{2}>0.$

\textbf{(ii) }Calculate the next iterate $x^{k+1}$ by%
\begin{equation}
x^{k+1}=\sum_{i=1}^{m}w_{i}y^{k+1,i}%
\end{equation}
where $w_{i}$ are fixed, user-chosen, positive weights with $\sum_{i=1}%
^{m}w_{i}=1.$
\end{algorithm}

The convergence analysis for this algorithm is currently available only for
consistent ($Q\neq\emptyset$) CFPs, see \cite{dossantos, moledo}. In our
experimental work, reported in the sequel, we applied Algorithm \ref{alg:ssp}
to CFPs without knowing whether or not they are consistent. Convergence is
diagnosed by performing plots of a \textit{proximity function }that measures
in some manner the infeasibility of the system. We used the weighted proximity
function of the form%
\begin{equation}
p(x):=(1/2)\sum_{i=1}^{m}w_{i}\parallel P_{i}(x)-x\parallel^{2}
\label{eq:proxfunct}%
\end{equation}
were $P_{i}(x)$ is the orthogonal projection of the point $x$ onto $Q_{i}.$ To
combat instabilities in those plots that appeared occasionally in our
experiments we used \textit{steering parameters} $\sigma_{k}$ instead of the
relaxation parameters $\alpha_{k}$ in Algorithm \ref{alg:ssp}$.$ To this end
we need the following definition.

\begin{definition}
A sequence $\{\sigma_{k}\}_{k=0}^{\infty}$ of real numbers $0\leq\sigma_{k}<1$
is called a \texttt{steering sequence} if it satisfies the following
conditions:%
\begin{equation}
\lim_{k\rightarrow\infty}\sigma_{k}=0, \label{steer1}%
\end{equation}%
\begin{equation}
\sum_{k=0}^{\infty}\sigma_{k}=+\infty, \label{steer2}%
\end{equation}%
\begin{equation}
\sum_{k=0}^{\infty}\mid\sigma_{k}-\sigma_{k+m}\mid<+\infty. \label{steer3}%
\end{equation}

\end{definition}

A historical and technical discussion of these conditions can be found in
\cite{bauschke96}. The sequential and simultaneous
Halpern-Lions-Wittmann-Bauschke (HLWB)\ algorithms discussed in Censor
\cite{hwlb-laa} employ the parameters of a steering sequence to ``force''
(steer) the iterates towards the solution of the best approximation problem
(BAP). This steering feature of the steering parameters has a profound effect
on the behavior of any sequence of iterates $\{x^{k}\}_{k=0}^{\infty}$. We
return to this point in Section \ref{sect:comp}.

\begin{algorithm}
\label{alg:steer}{The method of simultaneous subgradient projections }(SSP)
with steering{. }

\textbf{Initialization}: $x^{0}\in R^{n}$ is arbitrary.

\textbf{Iterative step}: \textbf{(i)} Given $x^{k},$ calculate, for all $i\in
I=\{1,2,\ldots,m\},$ intermediate iterates $y^{k+1,i}$ by
\begin{equation}
y^{k+1,i}=\left\{
\begin{array}
[c]{ll}%
x^{k}-\sigma_{k}\frac{\displaystyle f_{i}(x^{k})}{\displaystyle\parallel
\;t^{k}\parallel^{2}}t^{k}, & \mathrm{if}\;\;f_{i}(x^{k})>0,\\
x^{k}, & \mathrm{if}\;\;f_{i}(x^{k})\leq0,
\end{array}
\right.
\end{equation}
where $t^{k}\in\partial f_{i}(x^{k})$ is a subgradient of $f_{i}$ at the point
$x^{k}$, and $\{\sigma_{k}\}_{k=0}^{\infty}$ is a sequence of steering
parameters$.$

\textbf{(ii) }Calculate the next iterate $x^{k+1}$ by%
\begin{equation}
x^{k+1}=\sum_{i=1}^{m}w_{i}y^{k+1,i}%
\end{equation}
where $w_{i}$ are fixed, user-chosen, positive weights with $\sum_{i=1}%
^{m}w_{i}=1.$
\end{algorithm}

\section{Subgradient projections with strategical relaxation:
Preliminaries\label{sect:dan-prelim}}

Considering the CFP (\ref{(CFP problem)}), the \textit{envelope} of the family
of functions $\{f_{i}\}_{i=1}^{m}$ is the function%
\begin{equation}
f(x):=\max\{f_{i}(x)\mid i=1,2,\ldots,m\}
\end{equation}
which is also convex. Clearly, the consistent CFP is equivalent to finding a
point in%
\begin{equation}
Q=\cap_{i=1}^{m}Q_{i}=\{x\in R^{n}\mid f(x)\leq0\}.
\end{equation}

The subgradient projections algorithmic scheme that we propose here employs a
strategy for controlling the relaxation parameters in a specific manner,
leaving enough user-flexibility while giving some mathematical guarantees for
the algorithm's behavior in the inconsistent case. It is described as
follows:\medskip

\begin{algorithm}
\label{alg:dan}$\left.  {}\right.  $

\textbf{Initialization: }Let $M$ be a positive real number and let $x^{0}\in
R^{n}$ be any \textit{initial point}.

\textbf{Iterative step:} Given the current iterate $x^{k},$ set%
\begin{equation}
I(x^{k}):=\{i\mid1\leq i\leq m\text{ \ and \ }f_{i}(x^{k})=f(x^{k})\}
\label{Ix}%
\end{equation}
and choose a nonnegative vector $w^{k}=(w_{1}^{k},w_{2}^{k},\ldots,w_{m}%
^{k})\in R^{m}$ such that%
\begin{equation}
\sum_{i=1}^{m}w_{i}^{k}=1\text{ \ and \ }w_{i}^{k}=0\text{ \ if \ }i\notin
I(x^{k}).\text{\label{(W)}}%
\end{equation}
Let $\lambda_{k}$ be any nonnegative real number such that%
\begin{equation}
\max\left(  0,\text{ }f(x^{k})\right)  \leq\lambda_{k}M^{2}\leq2\max\left(
0,\text{ }f(x^{k})\right)  \label{lambda}%
\end{equation}
and calculate%
\begin{equation}
x^{k+1}=x^{k}-\lambda_{k}\sum_{i\in I(x^{k})}w_{i}^{k}\xi_{i}^{k},
\label{nextx}%
\end{equation}
where, for each $i\in I(x^{k}),$ we take a subgradient $\xi_{i}^{k}\in\partial
f_{i}(x^{k}).\medskip$
\end{algorithm}

It is interesting to note that any sequence $\left\{  x^{k}\right\}
_{k=0}^{\infty}$\ generated by this algorithm is well-defined, no matter how
$x^{0}$\ and $M$\ are chosen. Similarly to other algorithms described above,
Algorithm \ref{alg:dan} requires computing subgradients of convex functions.
In case a function is differentiable, this reduces to gradient calculations.
Otherwise, one can use the subgradient computing procedure presented in
Butnariu and Resmerita \cite{resmerita02}.

The procedure described above was previously studied in Butnariu and Mehrez
\cite{ButMeh}. The main result there shows that the procedure converges to a
solution of the CFP under two conditions: (i) that the solution set $Q$ has
nonempty interior and (ii) that the envelope $f$ is uniformly Lipschitz on
$R^{n},$ that is, there exists a positive real number $L$ such that%
\begin{equation}
\left\vert f(x)-f(y)\right\vert \leq L\left\Vert x-y\right\Vert ,\text{ \ for
all \ }x,y\in R^{n}. \label{GL}%
\end{equation}
Both conditions (i) and (ii) are restrictive and it is difficult to verify
their validity in practical applications. In the following we show that this
method converges to solutions of consistent CFPs under less demanding
conditions. In fact, we show that if the solution set $Q$ of the given CFP has
nonempty interior, then convergence of Algorithm \ref{alg:dan} to a point in
$Q$ is ensured even if the function $f$ is not uniformly Lipschitz on $R^{n}$
(i.e., even if $f$ does not satisfy condition (ii) above). However, verifying
whether $\operatorname{int}Q\neq\varnothing$ prior to solving a CFP may be
difficult or even impossible. Therefore, it is desirable to have alternative
conditions, which may be easier to verify in practice, that can ensure
convergence of our algorithm to solutions of the CFP, provided that such
solutions exist. This is why we prove convergence of Algorithm \ref{alg:dan}
to solutions of consistent CFPs whenever the envelope $f$ of the functions
$f_{i}$ involved in the given CFP is strictly convex. Strict convexity of the
envelope function $f$ associated with a consistent CFP implies that either the
solution set $Q$ of the CFP is a singleton, in which case $\operatorname{int}%
Q=\varnothing$, or that $Q$ contains (at least) two different solutions of the
CFP implying that int $Q\neq\varnothing$. Verification of whether $Q$ is a
singleton or not is as difficult as deciding whether int $Q\neq\varnothing$.
By contrast, since $f$ is strictly convex whenever each $f_{i}$ is strictly
convex, verification of strict convexity of $f$ may be relatively easily done
in some situations of practical interest, such as when each $f_{i}$ is a
quadratic convex function. In the latter case, strict convexity of $f_{i}$
amounts to positive definiteness of the matrix of its purely quadratic part.

It is interesting to note in this context that when the envelope $f$ of the
CFP is not strictly convex one may consider a \textquotedblleft
regularized\textquotedblright\ CFP in which each $f_{i}$ which is not strictly
convex is replaced by%
\begin{equation}
\overline{f_{i}}(x):=f_{i}(x)+\alpha\Vert x\Vert^{2}%
\end{equation}
for some positive real number $\alpha.$ Clearly, all $\overline{f_{i}}$ are
strictly convex and, thus, so is the envelope $\overline{f}$ of the
regularized problem. Therefore, if the regularized problem has solutions, then
our Algorithm \ref{alg:dan} will produce approximations of such solutions.
Moreover, any solution of the regularized problem is a solutions of the
original problem and, thus, by solving the regularized problem we implicitly
solve the original problem. The difficult part of this approach is that, even
if the original CFP is consistent then the regularized version of it may be
inconsistent for all, or for some, values $\alpha>0.$ How to decide whether an
$\alpha>0$ exists such that the corresponding regularized CFP is consistent
and how to compute such an $\alpha$ (if any) are questions whose answers we do
not know.

\section{Subgradient projections with strategical relaxation: convergence
analysis\label{sect:dan-converg}}

In order to discuss the convergence behavior of the subgradient projections
method with strategical relaxation, recall that convex functions defined on
the whole space $R^{n}$ are continuous and, consequently, are bounded on
bounded sets in $R^{n}$. Therefore, the application of Butnariu and Iusem
\cite[Proposition 1.1.11]{ButIusbook} or Bauschke and Borwein
\cite[Proposition 7.8]{bb96} to the convex function $f$ shows that it is
Lipschitz on bounded subsets of $R^{n},$ i.e., for any nonempty bounded subset
$S\subseteq R^{n}$ there exists a positive real number $L(S),$ called
\textit{a Lipschitz constant of }$f$\textit{\ over the set }$S,$ such that%
\begin{equation}
\left\vert f(x)-f(y)\right\vert \leq L(S)\left\Vert x-y\right\Vert ,\text{ for
all \ }x,y\in S. \label{LBS}%
\end{equation}
Our next result is a convergence theorem for Algorithm \ref{alg:dan} when
applied to a consistent CFP. It was noted in the previous section that
Algorithm \ref{alg:dan} is well-defined regardless of how the initial point
$x^{0}$ or the positive constant $M$ involved in the algorithm are chosen.
However, this is no guarantee that a sequences $\{x_{k}\}_{k=0}^{\infty}$
generated by Algorithm \ref{alg:dan} for random choices of $x^{0}$ and $M$
will converge to solutions of the CFP, even if such solutions exist. The
theorem below shows a way of choosing $x^{0}$ and $M$ which ensures that,
under some additional conditions for the problem data, the sequence
$\{x_{k}\}_{k=0}^{\infty}$ generated by Algorithm \ref{alg:dan} will
necessarily approximate a solution of the CFP (provided that solutions exist).
As shown in Section \ref{sec:implement} below, determining $x^{0}$ and $M$ as
required in the next theorem can be quite easily done for practically
significant classes of CFPs. Also, as shown in Section \ref{sect:comp},
determining $x^{0}$ and $M$ in this manner, enhances the self-adaptability of
the procedure to the problem data and makes Algorithm \ref{alg:dan} produce
approximations of solutions of the CFP which, in many cases, are more accurate
than those produced by other CFP solving algorithms.

\begin{theorem}
\label{thm:dan1}\textbf{ }\textit{If a positive number }$M$\textit{\ and an
initial point }$x^{0}$\textit{\ in Algorithm \ref{alg:dan} are chosen so that
}$M\geq L(B(x^{0},r))$\textit{\ for some positive real number }$r$%
\textit{\ satisfying the condition}
\begin{equation}
B(x^{0},r/2)\cap Q\neq\emptyset, \label{cond}%
\end{equation}
\textit{and if at least one of the following conditions holds:}

(i) $B(x^{0},r/2)\cap\operatorname*{int}\,Q\neq\emptyset,$

(ii) \textit{the function }$f$\textit{\ is strictly convex,}\newline%
\textit{then any sequence }$\left\{  x^{k}\right\}  _{k=0}^{\infty}%
,$\textit{\ generated by Algorithm \ref{alg:dan}, converges to an element of
}$Q.\smallskip$
\end{theorem}

We present the proof of Theorem \ref{thm:dan1} as a sequence of lemmas. To do
so, note that if $\left\{  x^{k}\right\}  _{k=0}^{\infty}$\textit{\ }is
generated by Algorithm \ref{alg:dan} then for each integer $k\geq0,$ we have%
\begin{equation}
x^{k+1}=x^{k}-\lambda_{k}\nu^{k},\label{xk}%
\end{equation}
where%
\begin{equation}
\nu^{k}:=\sum_{i\in I(x^{k})}w_{i}^{k}\xi_{i}^{k}\in\operatorname*{conv}%
\cup_{i\in I(x^{k})}\partial f_{i}(x^{k}).\label{niuk}%
\end{equation}
Using (\ref{xk}), for any $z\in R^{n},$ we have%
\begin{equation}
\left\Vert x^{k+1}-z\right\Vert ^{2}=\left\Vert x^{k}-z\right\Vert
^{2}+\lambda_{k}\left(  \lambda_{k}\left\Vert \nu^{k}\right\Vert
^{2}-2\left\langle \nu^{k},x^{k}-z\right\rangle \right)  .\label{basic}%
\end{equation}
By Clarke \cite[Proposition 2.3.12]{Clarke} we deduce that%
\begin{equation}
\partial f(x^{k})=\operatorname*{conv}\cup_{i\in I(x^{k})}\partial f_{i}%
(x^{k})\label{clarke}%
\end{equation}
and this implies that $\nu^{k}\in\partial f(x^{k})$ because of (\ref{niuk}).
Therefore,%
\begin{equation}
\left\langle \nu^{k},z-x^{k}\right\rangle \leq f_{+}^{\prime}(x^{k};z-x^{k})
\end{equation}
where $f_{+}^{\prime}(u;v)$ denotes the right-sided directional derivative at
$u$ in the direction $v.$ Now suppose that $M,$ $r$ and $x^{0}$ are chosen
according to the requirements of Theorem \ref{thm:dan1}, that is,%
\begin{equation}
r>0\text{, \ }M\geq L(B(x^{0},r))\text{ and }B(x^{0},r/2)\cap Q\neq
\emptyset.\label{prize}%
\end{equation}
Next we prove the following basic fact.

\begin{lemma}
\label{lem:dan1} \textit{If} (\ref{prize})\textit{\ is satisfied and if }$z\in
B(x^{0},r/2)\cap Q,$\textit{\ then for all} $k\geq0$, \textit{we have for any
sequence }$\left\{  x^{k}\right\}  _{k=0}^{\infty},$\textit{\ generated by
Algorithm \ref{alg:dan},}%
\begin{equation}
x^{k+1}\in B(x^{0},r)\text{ and }\left\|  x^{k+1}-z\right\|  \leq\left\|
x^{k}-z\right\|  \leq r/2\text{.} \label{bregmon}%
\end{equation}
\smallskip
\end{lemma}

\begin{proof}
We first show that if, for some integer $k\geq0$,%
\begin{equation}
x^{k}\in B(x^{0},r)\text{ and }\left\Vert x^{k}-z\right\Vert \leq r/2
\label{tnai}%
\end{equation}
then (\ref{bregmon}) holds. If $\lambda_{k}=0$ or $\nu^{k}=0,$ then, by
(\ref{xk}), we have $x^{k+1}=x^{k}$ which, combined with (\ref{tnai}), implies
(\ref{bregmon}). Assume now that $\lambda_{k}\neq0$ and $\nu^{k}\neq0.$ Since,
by (\ref{tnai}), $x^{k}\in B(x^{0},r),$ by (\ref{prize}) and by
\cite[Proposition 2.1.2(a)]{Clarke}, we deduce that%
\begin{equation}
M\geq L(B(x^{0},r))\geq\left\Vert \nu^{k}\right\Vert . \label{blabla}%
\end{equation}
According to (\ref{lambda}), we also have $f(x^{k})>0$ (otherwise $\lambda
_{k}=0$). Since $f(z)\leq0$ we obtain from the subgradient inequality%
\begin{equation}
\left\langle \nu^{k},x^{k}-z\right\rangle \geq f(x^{k})-f(z)\geq f(x^{k})>0.
\label{plus}%
\end{equation}
This and (\ref{blabla}) imply%
\begin{equation}
2\left\langle \nu^{k},x^{k}-z\right\rangle \geq2f(x^{k})\geq\lambda_{k}%
M^{2}\geq\lambda_{k}\left\Vert \nu^{k}\right\Vert ^{2}, \label{plusplus}%
\end{equation}
showing that the quantity inside the parentheses in (\ref{basic}) is
nonpositive. Thus, we deduce that%
\begin{equation}
\left\Vert x^{k+1}-z\right\Vert \leq\left\Vert x^{k}-z\right\Vert \leq r/2
\end{equation}
in this case too. This proves that if (\ref{tnai}) is true for all $k\geq0$,
then so is (\ref{bregmon}). Now, we prove by induction that (\ref{tnai}) is
true for all $k\geq0$. If $k=0$ then (\ref{tnai}) obviously holds. Suppose
that (\ref{tnai}) is satisfied for some $k=p.$ As shown above, this implies
that condition (\ref{bregmon}) is satisfied for $k=p$ and, thus, we have that%
\begin{equation}
x^{p+1}\in B(x^{0},r)\text{ and }\left\Vert x^{p+1}-z\right\Vert \leq r/2.
\end{equation}
Hence, condition (\ref{tnai}) also holds for $k=p+1.$ Consequently, condition
(\ref{bregmon}) holds for $k=p+1$ and this completes the proof.
\end{proof}

Observe that, according to Lemma \ref{lem:dan1}, if $\left\{  x^{k}\right\}
_{k=0}^{\infty}$ is a sequence generated by Algorithm \ref{alg:dan} and if the
conditions (\ref{prize}) are satisfied, then there exists $z\in B(x^{0}%
,r/2)\cap Q$ and for any such $z$ the sequence $\left\{  \left\Vert
x^{k}-z\right\Vert \right\}  _{k=0}^{\infty}$ is nonincreasing and bounded
from below and, therefore, convergent. Since the sequence $\left\{  \left\Vert
x^{k}-z\right\Vert \right\}  _{k=0}^{\infty}$ is convergent it is also bounded
and, consequently, the sequence $\left\{  x^{k}\right\}  _{k=0}^{\infty}$ is
bounded too. This shows that the next result applies to any sequence $\left\{
x^{k}\right\}  _{k=0}^{\infty}$ generated by Algorithm \ref{alg:dan} under the
assumptions of Theorem \ref{thm:dan1}.

\begin{lemma}
\label{lem:subseq} I\textit{f }$\left\{  x^{k}\right\}  _{k=0}^{\infty}%
$\textit{\ is a bounded sequence generated by Algorithm \ref{alg:dan}, then
the sequence }$\left\{  x^{k}\right\}  _{k=0}^{\infty}$\textit{\ has
accumulation points and for each accumulation point }$x^{\ast}$\textit{\ of
}$\left\{  x^{k}\right\}  _{k=0}^{\infty}$\textit{\ there exists a sequence of
natural numbers }$\left\{  k_{s}\right\}  _{s=0}^{\infty}$\textit{\ such that
the following limits exist}%
\begin{align}
x^{\ast} &  =\lim_{s\rightarrow\infty}x^{k_{s}},\text{ }\lambda_{\ast}%
=\lim_{s\rightarrow\infty}\lambda_{k_{s}},\label{pips}\\
\xi_{i}^{\ast} &  =\lim_{s\rightarrow\infty}\xi_{i}^{k_{s}},\text{ }%
w_{i}^{\ast}=\lim_{s\rightarrow\infty}w_{i}^{k_{s}},\text{ for all
}i=1,2,\ldots,m,\\
\nu^{\ast} &  =\lim_{s\rightarrow\infty}\nu^{k_{s}},\label{eq:pips3}%
\end{align}
\textit{and we have}%
\begin{equation}
w^{\ast}:=(w_{1}^{\ast},w_{2}^{\ast},\ldots,w_{m}^{\ast})\in R_{+}%
^{m}\ \text{and }\sum_{i\in I(x^{\ast})}w_{i}^{\ast}=1\label{legaturi1}%
\end{equation}
and%
\begin{equation}
\nu^{\ast}=\sum_{i\in I(x^{\ast})}w_{i}^{\ast}\xi_{i}^{\ast}\in\partial
f(x^{\ast}).\label{legaturi2}%
\end{equation}
\textit{Moreover, if }$\lambda_{\ast}=0,$\textit{\ then }$x^{\ast}%
$\textit{\ is a solution of the CFP}.\smallskip
\end{lemma}

\begin{proof}
The sequence $\left\{  x^{k}\right\}  _{k=0}^{\infty}$ is bounded and, thus,
has accumulation points.\ Let $x^{\ast}$ be an accumulation point of $\left\{
x^{k}\right\}  _{k=0}^{\infty}$ and let $\left\{  x^{p_{s}}\right\}
_{s=0}^{\infty}$ be a convergent subsequence of $\left\{  x^{k}\right\}
_{k=0}^{\infty}$ such that $x^{\ast}=\lim_{s\rightarrow\infty}x^{p_{s}}.$ The
function $f$ is continuous (since it is real-valued and convex on $R^{n}$),
hence, it is bounded on bounded subsets of $R^{n}$. Therefore, the sequence
$\left\{  f(x^{p_{s}})\right\}  _{s=0}^{\infty}$ converges to $f(x^{\ast})$
and the sequence $\left\{  f(x^{k})\right\}  _{k=0}^{\infty}$ is bounded. By
(\ref{lambda}), boundedness of $\left\{  f(x^{k})\right\}  _{k=0}^{\infty}$
implies that the sequence $\left\{  \lambda_{k}\right\}  _{k=0}^{\infty}$ is
bounded. Since, for every $i=1,2,\ldots,m,$ the operator $\partial f_{i}%
:R^{n}\rightarrow2^{R^{n}}$ is monotone, it is locally bounded (cf. Pascali
and Sburlan \cite[Theorem on p. 104]{PasSbu}).

Consequently, there exists a neighborhood $U$ of $x^{\ast}$ on which all
$\partial f_{i},$ $i=1,2,\ldots,m,$ are bounded. Clearly, since $x^{\ast}%
=\lim_{s\rightarrow\infty}x^{p_{s}}$, the neighborhood $U$ contains all but
finitely many terms of the sequence $\left\{  x^{p_{s}}\right\}
_{s=0}^{\infty}.$ This implies that the sequences $\left\{  \xi_{i}^{p_{s}%
}\right\}  _{s=0}^{\infty}$ are uniformly bounded and, therefore, the sequence
$\left\{  \nu^{p_{s}}\right\}  _{s=0}^{\infty}$ is bounded too.

Therefore, there exist a subsequence $\left\{  k_{s}\right\}  _{s=0}^{\infty}$
of $\left\{  p_{s}\right\}  _{s=0}^{\infty}$ such that the limits in
(\ref{pips})--(\ref{eq:pips3}) exist. Obviously, the vector $w^{\ast}%
=(w_{1}^{\ast},w_{2}^{\ast},\ldots,w_{m}^{\ast})\in R_{+}^{m},$ and, according
to \cite[Lemma 1]{ButMeh}, we also have $\sum_{i\in I(x^{\ast})}w_{i}^{\ast
}=1$. This and (\ref{niuk}) imply that $\nu^{\ast}=\sum_{i\in I(x^{\ast}%
)}w_{i}^{\ast}\xi_{i}^{\ast}.$

Observe that, since $\nu^{k_{s}}\in\partial f(x^{k_{s}})$ for all $s\geq0$,
and since $\partial f$ is a closed mapping (cf. Phelps \cite[Proposition
2.5]{Phe}), we have that $\nu^{\ast}\in\partial f(x^{\ast}).$ Now, if
$\lambda_{\ast}=0$ then, according to (\ref{lambda}), and the continuity of
$f,$ we deduce%
\begin{equation}
0\leq\max\{0,f(x^{\ast})\}=\lim_{s\rightarrow\infty}\max\{0,f(x^{k_{s}}%
)\}\leq\lim_{s\rightarrow\infty}\lambda_{k_{s}}M^{2}=\lambda_{\ast}M^{2}=0,
\end{equation}
which implies that $f(x^{\ast})\leq0$, that is, $x^{\ast}\in Q.$
\end{proof}

\begin{lemma}
\label{lem:5}\textit{Let }$\left\{  x^{k}\right\}  _{k=0}^{\infty}%
$\textit{\ be a sequence generated by Algorithm \ref{alg:dan}. If
(\ref{prize}) is satisfied and if at least one of the conditions (i) or (ii)
of Theorem \ref{thm:dan1} holds, then the sequence }$\left\{  x^{k}\right\}
_{k=0}^{\infty}$\textit{\ has accumulation points and any such point belongs
to }$Q.\smallskip$
\end{lemma}

\begin{proof}
As noted above, when (\ref{prize})\ is satisfied then the sequence $\left\{
x^{k}\right\}  _{k=0}^{\infty}$ is bounded and, hence, it has accumulation
points. Let $x^{\ast}$ be such an accumulation point and let $\left\{
k_{s}\right\}  _{s=0}^{\infty}$ be the sequence of natural numbers associated
with $x^{\ast}$ whose existence is guaranteed by Lemma \ref{lem:subseq}.
Since, for any $z\in C\cap B(x^{0},r/2),$ the sequence $\left\{  \left\Vert
x^{k}-z\right\Vert \right\}  _{k=0}^{\infty}$ is convergent (cf. Lemma
\ref{lem:dan1}) we deduce that%
\begin{align}
\left\Vert x^{\ast}-z\right\Vert  &  =\lim_{s\rightarrow\infty}\left\Vert
x^{k_{s}}-z\right\Vert =\lim_{k\rightarrow\infty}\left\Vert x^{k}-z\right\Vert
=\lim_{s\rightarrow\infty}\left\Vert x^{k_{s}+1}-z\right\Vert \\
&  =\left\Vert x^{\ast}-\lambda_{\ast}\nu^{\ast}-z\right\Vert .\nonumber
\end{align}
This implies%
\begin{equation}
\left\Vert x^{\ast}-z\right\Vert ^{2}=\left\Vert x^{\ast}-z\right\Vert
^{2}+\lambda_{\ast}\left(  \lambda_{\ast}\left\Vert \nu^{\ast}\right\Vert
^{2}-2\left\langle \nu^{\ast},x^{\ast}-z\right\rangle \right)  .\label{kill}%
\end{equation}
If $\lambda_{\ast}=0,$ then $x^{\ast}\in Q$ by Lemma \ref{lem:subseq}. Suppose
that $\lambda_{\ast}>0.$ Then, by (\ref{kill}), we have
\begin{equation}
\lambda_{\ast}\left\Vert \nu^{\ast}\right\Vert ^{2}-2\left\langle \nu^{\ast
},x^{\ast}-z\right\rangle =0,\label{protea}%
\end{equation}
for all $z\in C\cap B(x^{0},r/2).$ We distinguish now between two possible cases.

\textbf{Case I}: Assume that condition (\textit{i}) of Theorem \ref{thm:dan1}
is satisfied. According to (\ref{protea}), the set $Q\cap B(x^{0},r/2)$ is
contained in the hyperplane%
\begin{equation}
H:=\left\{  x\in R^{n}\mid\left\langle \nu^{\ast},x\right\rangle =(1/2)\left(
2\left\langle \nu^{\ast},x^{\ast}\right\rangle -\lambda_{\ast}\left\|
\nu^{\ast}\right\|  ^{2}\right)  \right\}  .
\end{equation}
By \textit{condition (i) of Theorem \ref{thm:dan1}}, it follows that
$\mathrm{\operatorname*{int}}\,\left(  Q\cap B(x^{0},r/2)\right)
\neq\emptyset$ and this is an open set contained in
$\mathrm{\operatorname*{int}}\,H.$ So, unless $\nu^{\ast}=0$ (in which case
$H=R^{n}$), we have reached a contradiction because
$\mathrm{\operatorname*{int}}\,H=\emptyset.$ Therefore, we must have
$\nu^{\ast}=0.$ According to Lemma \ref{lem:subseq}, we have $0=\nu^{\ast}%
\in\partial f(x^{\ast})$ which implies that $x^{\ast}$ is a global minimizer
of $f.$ Consequently, for any $z\in Q$ we have $f(x^{\ast})\leq f(z)\leq0,$
that is, $x^{\ast}\in Q.$

\textbf{Case II}: Assume that condition (\textit{ii}) of Theorem
\ref{thm:dan1} is satisfied. According to (\ref{protea}), we have%
\begin{equation}
\lambda_{\ast}\left\Vert \nu^{\ast}\right\Vert ^{2}=2\left\langle \nu^{\ast
},x^{\ast}-z\right\rangle . \label{muc}%
\end{equation}
By (\ref{lambda}), the definition of $M$ and \cite[Proposition 2.1.2]{Clarke}
we deduce that%
\begin{equation}
2f(x^{k_{s}})\geq\lambda_{k_{s}}M^{2}\geq\lambda_{k_{s}}\left\Vert \nu^{k_{s}%
}\right\Vert ^{2},
\end{equation}
for all integers $s\geq0$. Letting $s\rightarrow\infty$ we get%
\begin{equation}
2f(x^{\ast})\geq\lambda_{\ast}M^{2}\geq\lambda_{\ast}\left\Vert \nu^{\ast
}\right\Vert ^{2}=2\left\langle \nu^{\ast},x^{\ast}-z\right\rangle ,
\end{equation}
where the last equality follows from (\ref{muc}). Consequently, we have%
\begin{equation}
f(x^{\ast})\geq\left\langle \nu^{\ast},x^{\ast}-z\right\rangle ,\text{ \ \ for
all \ \ }z\in Q\cap B(x^{0},r/2).
\end{equation}
Convexity of $f$ implies that, for all $z\in Q\cap B(x^{0},r/2),$%
\begin{equation}
-f(x^{\ast})\leq\left\langle \nu^{\ast},z-x^{\ast}\right\rangle \leq
f(z)-f(x^{\ast})\leq-f(x^{\ast}).
\end{equation}
Therefore, we have that%
\begin{equation}
-f(x^{\ast})=\left\langle \nu^{\ast},z-x^{\ast}\right\rangle =f(z)-f(x^{\ast
}),\text{ for all \ \ }z\in Q\cap B(x^{0},r/2).
\end{equation}
Thus $f(z)=0,$ for all $z\in Q\cap B(x^{0},r/2).$ Hence, using again the
convexity of $f,$ we deduce that, for all \ \ $z\in Q\cap B(x^{0},r/2),$%
\begin{equation}
f_{+}^{\prime}(x^{\ast};z-x^{\ast})\leq f(z)-f(x^{\ast})=-f(x^{\ast
})=\left\langle \nu^{\ast},z-x^{\ast}\right\rangle \leq f_{+}^{\prime}%
(x^{\ast};z-x^{\ast}).
\end{equation}
This implies%
\begin{equation}
f_{+}^{\prime}(x^{\ast};z-x^{\ast})=\left\langle \nu^{\ast},z-x^{\ast
}\right\rangle =f(z)-f(x^{\ast}),\text{ for all \ \ }z\in Q\cap B(x^{0},r/2).
\label{mule}%
\end{equation}
Since, by \textit{condition (ii) of Theorem \ref{thm:dan1}}, $f$ is strictly
convex, we also have (see \cite[Proposition 1.1.4]{ButIusbook}) that%
\begin{equation}
f_{+}^{\prime}(x^{\ast};z-x^{\ast})<f(z)-f(x^{\ast}),\text{ for all \ \ }%
z\in\left(  Q\cap B(x^{0},r/2)\right)  \backslash\{x^{\ast}\}.
\end{equation}
Hence, the equalities in (\ref{mule}) cannot hold unless $Q\cap B(x^{0}%
,r/2)=\{x^{\ast}\}$ and, thus, $x^{\ast}\in Q.$
\end{proof}

The previous lemmas show that if (\ref{prize}) holds and if one of the
conditions \textit{(i)} or \textit{(ii)} of Theorem \ref{thm:dan1}\textit{ }is
satisfied, then the sequence $\left\{  x^{k}\right\}  _{k=0}^{\infty}$ is
bounded and all its accumulation points are in $Q.$ In fact, the results above
say something more. Namely, in view of Lemma \ref{lem:dan1}, they show that if
(\ref{prize}) holds and if one of the conditions \textit{(i)} or \textit{(ii)}
of Theorem \ref{thm:dan1} is satisfied, then all accumulation points $x^{\ast
}$ of $\left\{  x^{k}\right\}  _{k=0}^{\infty}$ are contained in $Q\cap
B(x^{0},r)$ because all $x^{k}$ are in $B(x^{0},r)$ by (\ref{bregmon}). In
order to complete the proof of Theorem \ref{thm:dan1}, it remains to show that
the following result is true.\smallskip

\begin{lemma}
\textit{Under the conditions} of Theorem \ref{thm:dan1} \textit{any sequence
}$\left\{  x^{k}\right\}  _{k=0}^{\infty}$,\ generated by Algorithm
\ref{alg:dan},\textit{\ has at most one accumulation point.\smallskip}
\end{lemma}

\begin{proof}
Observe that, under the conditions of Theorem \ref{thm:dan1} the conditions
(\ref{prize}) are satisfied and, therefore, the sequence $\left\{
x^{k}\right\}  _{k=0}^{\infty}$ is bounded. Let $x^{\ast}$ be an accumulation
point of $\left\{  x^{k}\right\}  _{k=0}^{\infty}.$ By Lemma \ref{lem:5} we
deduce that $x^{\ast}\in Q,$ i.e., $f(x^{\ast})\leq0.$ Consequently, for any
natural number $k$ we have%
\[
\left\langle \nu^{k},x^{k}-x^{\ast}\right\rangle \geq f(x^{k})-f(x^{\ast})\geq
f(x^{k}).
\]
Now, using this fact, a reasoning similar to that which proves (\ref{plusplus}%
) but made with $x^{\ast}$ instead of $z$ leads to
\[
2\left\langle \nu^{k},x^{k}-x^{\ast}\right\rangle \geq\lambda_{k}\left\Vert
\nu^{k}\right\Vert ^{2},
\]
for all natural numbers $k.$ This and (\ref{basic}) combined imply that the
sequence $\left\{  \left\Vert x^{k}-x^{\ast}\right\Vert \right\}
_{k=0}^{\infty}$ is non-increasing and, therefore, convergent. Consequently,
if $\left\{  x^{k_{p}}\right\}  _{p=0}^{\infty}$ is a subsequence of $\left\{
x^{k}\right\}  _{k=0}^{\infty}$ such that $\lim_{p\rightarrow\infty}x^{k_{p}%
}=x^{\ast},$ we have%
\[
\lim_{k\rightarrow\infty}\left\Vert x^{k}-x^{\ast}\right\Vert =\lim
_{p\rightarrow\infty}\left\Vert x^{k_{p}}-x^{\ast}\right\Vert =0,
\]
showing that any accumulation point $x^{\ast}$ of $\left\{  x^{k}\right\}
_{k=0}^{\infty}$ is exactly the limit of $\left\{  x^{k}\right\}
_{k=0}^{\infty}.$
\end{proof}

The application of Theorem \ref{thm:dan1} depends on our ability to choose
numbers $M$ and $r$ and a vector $x^{0}$ such that condition (\ref{prize}) is
satisfied. We show below that this can be done when the functions $f_{i}$ of
the CFP (\ref{(CFP problem)}) are quadratic or affine and there is some a
priori known ball which intersects $Q$. In actual applications it may be
difficult to a priori decide whether the CFP (\ref{(CFP problem)}) has or does
not have solutions. However, as noted above, Algorithm \ref{alg:dan} is
well-defined and will generate sequences $\left\{  x^{k}\right\}
_{k=0}^{\infty}$ no matter how the initial data $M,$ $r$ and $x^{0}$ are
chosen. This leads to the question whether it is possible to decide if $Q$ is
empty or not by simply analyzing the behavior of sequences $\left\{
x^{k}\right\}  _{k=0}^{\infty}$ generated by Algorithm \ref{alg:dan}. A
partial answer to this question is contained in the following result.

\begin{corollary}
\label{cor:7}\textit{Suppose that the CFP (\ref{(CFP problem)}) has no
solution and that the envelope }$f$\textit{\ is strictly convex. Then, no
matter how the initial vector }$x^{0}$\textit{\ and the positive number }%
$M$\textit{\ are chosen, any sequence }$\left\{  x^{k}\right\}  _{k=0}%
^{\infty},$\textit{\ generated by Algorithm \ref{alg:dan}, has the following
properties:}

(i) \textit{If }$\left\{  x^{k}\right\}  _{k=0}^{\infty}$\textit{\ is bounded
and}%
\begin{equation}
\lim_{k\rightarrow\infty}\left\|  x^{k+1}-x^{k}\right\|  =0, \label{limit}%
\end{equation}
\textit{then }$f$ \textit{has a (necessarily unique) minimizer and} $\left\{
x^{k}\right\}  _{k=0}^{\infty}$ \textit{converges to that minimizer while}%
\begin{equation}
\lim_{k\rightarrow\infty}f(x^{k})=\inf\{f(x)\mid x\in R^{n}\}. \label{val}%
\end{equation}

(ii) \textit{If }$f$\textit{\ has no minimizer then the sequence }$\left\{
x^{k}\right\}  _{k=0}^{\infty}$\textit{\ is unbounded or the sequence
}$\left\{  \left\Vert x^{k+1}-x^{k}\right\Vert \right\}  _{k=0}^{\infty}$
\textit{does not converge to zero.}\smallskip
\end{corollary}

\begin{proof}
Clearly, (\textit{ii}) is a consequence of (\textit{i}). In order to prove
(\textit{i}) observe that, since the CFP (\ref{(CFP problem)}) has no
solution, all values of $f$ are positive. Also, if $f$ has a minimizer, then
this minimizer is unique because $f$ is strictly convex.

If $\left\{  x^{k}\right\}  _{k=0}^{\infty}$ is bounded then it has an
accumulation point, say, $x^{\ast}$. By Lemma \ref{lem:subseq} there exists a
sequence of positive integers $\left\{  k_{s}\right\}  _{s=0}^{\infty}$ such
that (\ref{pips}) and (\ref{legaturi1})--(\ref{legaturi2}) are satisfied.
Using Lemma \ref{lem:subseq} again, we deduce that, if the limit
$\lambda_{\ast}$ in (\ref{pips}) is zero, then the vector $x^{\ast}%
=\lim_{s\rightarrow\infty}x^{k_{s}}$ is a solution of the CFP
(\ref{(CFP problem)}), i.e., $f(x^{\ast})\leq0,$ contradicting the assumption
that the CFP (\ref{(CFP problem)}) has no solution. Hence, $\lambda_{\ast}>0.$
By (\ref{xk}), (\ref{limit}) and (\ref{pips}) we have that%
\begin{equation}
0=\lim_{s\rightarrow\infty}\lambda_{k_{s}}\nu^{k_{s}}=\lambda_{\ast}\nu^{\ast
}.
\end{equation}
Thus, we deduce that $\nu^{\ast}=0.$ From (\ref{legaturi1})--(\ref{legaturi2})
and \cite[Proposition 2.3.12]{Clarke} we obtain%
\begin{equation}
0=\nu^{\ast}=\sum_{i\in I(x^{\ast})}w_{i}^{\ast}\xi_{i}^{\ast}\in\partial
f(x^{\ast}),
\end{equation}
showing that $x^{\ast}$ is a minimizer of $f.$ So, all accumulation points of
$\left\{  x^{k}\right\}  _{k=0}^{\infty}$ coincide because $f$ has no more
than one minimizer. Consequently, the bounded sequence $\left\{
x^{k}\right\}  _{k=0}^{\infty}$ converges and its limit is the unique
minimizer of $f.$
\end{proof}

\begin{remark}
Checking numerically a condition such as (ii) in Corollary \ref{cor:7} or the
condition in Corollary \ref{cor:10} below seems virtually impossible. But
there is no escape from such situations in such mathematically-oriented
results. Condition (ii) in Corollary \ref{cor:7} is meaningful in the
inconsistent case in which a feasible point does not exist but a proximity
function that \textquotedblleft measures\textquotedblright\ the feasibility
violation of the limit point can be minimized. An easy adaptation of the proof
of Corollary 4.6 shows that if the sequence $\left\{  x^{k}\right\}
_{k=0}^{\infty}$ has a bounded subsequence $\left\{  x^{k_{t}}\right\}
_{t=0}^{\infty}$ such that the limit $\lim_{t\rightarrow\infty}(x^{k_{t}%
+1}-x^{k_{t}})=0$, then all accumulation points of $\left\{  x^{k_{t}%
}\right\}  _{t=0}^{\infty}$ are minimizers of $f$ (even if $f$ happens to be
not strictly convex).
\end{remark}

\begin{remark}
The fact that for some choice of $x^{0}$ and $M$ a sequence $\left\{
x^{k}\right\}  _{k=0}^{\infty},$ generated by \textit{Algorithm \ref{alg:dan}%
,} has the property that $\lim_{k\rightarrow\infty}f(x^{k})=0,$ does not imply
that the CFP (\ref{(CFP problem)}) has a solution. For example, take in
(\ref{(CFP problem)}) $m=n=1$ and $f_{1}(x)=e^{-x}.$ Clearly, in this case
(\ref{(CFP problem)}) has no solution and $f=f_{1}$. However, for $x^{0}=0$,
$M=1$ and $\lambda_{k}=(3/2)f(x^{k}),$ we have $\lim_{k\rightarrow\infty
}f(x^{k})=0.$
\end{remark}

A meaningful implication of Corollary \ref{cor:7} is the following result.

\begin{corollary}
\label{cor:10}Suppose that the CFP (2) has no solution and that $f$ is
strictly convex. Then, no matter how the initial vector $x^{0}$ and the
positive number $M$ are chosen in Algorithm \ref{alg:dan}, the following
holds: If the series $\sum_{k=0}^{\infty}\left\|  x^{k}-x^{k+1}\right\|  $
converges, then the function $f$ has a unique global minimizer and the
sequence $\left\{  x^{k}\right\}  _{k=0}^{\infty},$ generated by Algorithm
\ref{alg:dan}, converges to that minimizer while the sequence $\left\{
f(x^{k})\right\}  _{k=0}^{\infty}$ converges to $\inf\{f(x)\mid x\in R^{n}\}.$
\end{corollary}

\begin{proof}
When $\sum_{k=0}^{\infty}\left\Vert x^{k}-x^{k+1}\right\Vert $ converges to
some number $S$ we have%
\begin{equation}
\left\Vert x^{0}-x^{k+1}\right\Vert \leq\sum_{\ell=0}^{k}\left\Vert x^{\ell
}-x^{\ell+1}\right\Vert \leq S,
\end{equation}
for all integers $k\geq0.$ This implies that the sequence $\left\{
x^{k}\right\}  _{k=0}^{\infty}$ is bounded and $\lim_{k\rightarrow\infty
}\left\Vert x^{k}-x^{k+1}\right\Vert =0$. Hence, by applying Corollary 4.6, we
complete the proof.
\end{proof}

\begin{remark}
\label{rem:a}Finding an initial vector $x^{0},$ the radius $r$ and a positive
number $M$ satisfying condition $M\geq L(B(x^{0},r))$ (and satisfying
(\ref{cond}) provided that $Q$ is nonempty) when there is no a priori
knowledge about \ the existence of a solution of the CFP can be quite easily
done when at least one of the sets $Q_{i},$ say $Q_{i_{0}},$ is bounded and
the functions $f_{i}$ are differentiable. In this case it is sufficient to
determine a vector $x^{0}$ and a positive number $r$ large enough so that the
ball $B(x^{0},r/2)$ contains $Q_{i_{0}}.$ Clearly, for such a ball, if $Q$ is
nonempty, then condition (\ref{cond}) holds. Once the ball $B(x^{0},r)$ is
determined, finding a number $M\geq L(B(x^{0},r))$ can be done by taking into
account that the gradients of the differentiable convex functions $f_{i}%
:R^{n}\rightarrow R$ are necessarily continuous and, therefore, the numbers%
\begin{equation}
L_{i}=\sup\{\left\Vert \nabla f_{i}(x)\right\Vert \mid x\in B(x^{0},r)\}
\end{equation}
are necessarily finite. Since $L:=\max\{L_{i}\mid1\leq i\leq m\}$ is
necessarily a Lipschitz constant of $f$ over $B(x^{0},r),$ one can take $M=L$.
\end{remark}

\begin{remark}
\label{rem:b}The method of choosing $x^{0},$ $r$ and $M$ presented in Remark
\ref{rem:a} does not require a priori knowledge of the existence of a solution
of the CFP and can be applied even when $Q$ is empty. In such a case one
should compute, along the iterative procedure of Algorithm \ref{alg:dan}, the
sums $S_{k}=\sum_{\ell=0}^{k}\left\Vert x^{\ell}-x^{\ell+1}\right\Vert .$
Theorem \ref{thm:dan1} and Corollary \ref{cor:10} then provide the following
insights and tools for solving the CFP, provided that $f$ is strictly convex:

\begin{itemize}
\item If along the computational process the sequence $S_{k}$ remains bounded
from above by some number $S^{\ast}$ while the sequence $\{f(x^{k}%
)\}_{k=0}^{\infty}$ stabilizes itself asymptotically at some \textbf{positive}
value, then the given CFP has no solution, but the sequence $\left\{
x^{k}\right\}  _{k=0}^{\infty}$ still approximates a global minimum of $f$
which may be taken as a surrogate solution of the given CFP.

\item If along the computational process the sequence $S_{k}$ remains bounded
from above by some number $S^{\ast}$ while the sequence $\{f(x^{k}%
)\}_{k=0}^{\infty}$ stabilizes itself asymptotically at some
\textbf{nonpositive} value, then the given CFP has a solution, and the
sequence $\left\{  x^{k}\right\}  _{k=0}^{\infty}$ approximates such a
solution. 
\end{itemize}
\end{remark}

\section{Implementation of Algorithm \ref{alg:dan} for linear or quadratic
functions\label{sec:implement}}

Application of Algorithm \ref{alg:dan} does not require knowledge of the
constant $r.$ However, in order to implement this algorithm, such that the
conditions for convergence will be guaranteed, we have to determine numbers
$r$ and $M$ required by Theorem \ref{thm:dan1}. The method proposed in Remark
\ref{rem:a} might yield a very large value of $r.$ This is due to the
mathematical generality of Remark \ref{rem:a}. The quadratic and affine cases
treated next seem to be restrictive from the theoretical/mathematical point of
view, but their importance lies in the fact that they cover many significant
real-world applications.

We deal first with the problem of determining a number $M$ such that%
\begin{equation}
M\geq L(B(x^{0},r)), \label{kiki}%
\end{equation}
provided that an $r>0$ is given. Recall that if $g:R^{n}\rightarrow R$ is a
continuously differentiable function then, by Taylor's formula, we have that,
whenever $x,y\in B(x^{0},r),$ there exists a $u\in\left[  x,y\right]  $ such
that%
\begin{align}
&  \mid g(y)-g(x)\mid=\mid\left\langle \nabla g(u),y-x\right\rangle \mid
\leq\left\Vert \nabla g(u)\right\Vert \left\Vert y-x\right\Vert \nonumber\\
&  \leq\left\Vert y-x\right\Vert \max\{\left\Vert \nabla g(u)\right\Vert \mid
u\in B(x^{0},r)\}. \label{g}%
\end{align}
This shows that%
\begin{equation}
\max\{\left\Vert \nabla g(u)\right\Vert \mid u\in B(x^{0},r)\} \label{eq:Lip}%
\end{equation}
is a Lipschitz constant for $g$ on $B(x^{0},r)$. Suppose now that each
function $f_{i}$ is either linear or quadratic. Denote $I_{1}=\{i\mid1\leq
i\leq m,$ $\ f_{i}$ is linear$\}$ and $I_{2}=\{i\mid1\leq i\leq m,$ $\ f_{i}$
is quadratic$\}.$ Namely,%
\begin{equation}
f_{i}(x)=\left\langle a^{i},x\right\rangle +b_{i},\text{ \ for all \ }i\in
I_{1},
\end{equation}
with $a^{i}\in R^{n}\backslash\left\{  0\right\}  $ and $b_{i}\in R$, and%
\begin{equation}
f_{i}(x)=\left\langle x,U_{i}x\right\rangle +\left\langle a^{i},x\right\rangle
+b_{i},\text{ for all \ }i\in I_{2},
\end{equation}
where $U_{i}=(u_{\ell,k}^{i})$ is a symmetric positive semidefinite $n\times
n$ matrix, $a^{i}\in R^{n}$ and $b_{i}\in R$. We have, of course,%
\begin{equation}
\nabla f_{i}(x)=\left\{
\begin{array}
[c]{cc}%
a^{i}, & \text{if \ }i\in I_{1},\\
2U_{i}x+a^{i}, & \text{if \ }i\in I_{2},
\end{array}
\right.
\end{equation}
so that (\ref{eq:Lip}) can give us Lipschitz constants for each $f_{i}$ over
$B(x^{0},r).$ Denote%
\begin{equation}
L_{i}:=\left\{
\begin{array}
[c]{cc}%
\left\Vert a^{i}\right\Vert , & \text{if }i\in I_{1},\\
2\left\Vert U_{i}\right\Vert _{\infty}\left(  \left\Vert x^{0}\right\Vert
+r\right)  +\left\Vert a^{i}\right\Vert , & \text{if }i\in I_{2},
\end{array}
\right.  \label{eq:Li}%
\end{equation}
where $\left\Vert U_{i}\right\Vert _{\infty}$ is the operator norm of $U_{i}.$
Due to (\ref{clarke}), this implies that $\cup_{x\in B(x^{0},r)}\partial
f(x)\subseteq B(0,L)$ where%
\begin{equation}
L:=\max\{L_{i}\mid1\leq i\leq m\}.
\end{equation}
Taking $\xi\in\partial f(x)$ and $\zeta\in\partial f(y),$ for some $x,y\in
B(x^{0},r)$, we have%
\begin{align}
L\left\Vert x-y\right\Vert  &  \geq\left\Vert \zeta\right\Vert \left\Vert
x-y\right\Vert \geq\left\langle \zeta,y-x\right\rangle \geq
f(y)-f(x)\nonumber\\
&  \geq\left\langle \xi,y-x\right\rangle \geq-\left\Vert \xi\right\Vert
\left\Vert x-y\right\Vert \geq-L\left\Vert x-y\right\Vert ,
\end{align}
which implies%
\begin{equation}
\left\vert f(y)-f(x)\right\vert \leq L\left\Vert x-y\right\Vert ,\text{ for
all \ }x,y\in B(x^{0},r).
\end{equation}
In other words, $L$ is a Lipschitz constant of $f$ over $B(x^{0},r).$ Thus,
given an $r>0,$ we can take $M$ to be any number such that $M\geq L.$ Note
that choosing $x^{0}$ such that the corresponding $r$ is small may speed up
the computational process by reducing the number of iterations needed to reach
a reasonably good approximate solution of the CFP. In general, determining a
number $r$ is straightforward when one has some information about the range of
variation of the coordinates of some solutions to the CFP.

For instance, if one knows a priori that the solutions of the CFP are vectors
$x=(x_{j})_{j=1}^{n}$ such that%
\begin{equation}
\ell_{j}\leq x_{j}\leq u_{j},\text{ }1\leq j\leq n, \label{neweq1}%
\end{equation}
where, $\ell_{j},u_{j}\in R$, for all $j,$ then the set $Q$ is contained in
the hypercube of edge length $\delta=u_{\max}-\ell_{\min},$ whose faces are
parallel to the axes of the coordinates, and centered at the point $x^{0}$
whose coordinates are $x_{j}^{0}=\textstyle%
\frac{1}{2}%
(\ell_{\min}+u_{\max}),$ where%
\begin{equation}
\ell_{\min}:=\min\{\ell_{j}\mid1\leq j\leq n\}\text{ and }u_{\max}%
:=\max\{u_{j}\mid1\leq j\leq n\}.
\end{equation}
Therefore, by choosing this $x^{0}$ as the initial point for Algorithm
\ref{alg:dan} and choosing $r=\sqrt{n}\delta,$ condition (\ref{cond}) holds.

\section{Computational results\label{sect:comp}}

In this section, we compare the performance of Algorithms \ref{alg:ssp},
\ref{alg:steer} and \ref{alg:dan} by examining a few test problems. There are
a number of degrees-of-freedom used to evaluate and compare the performance of
the algorithms. These are the maximum number of iterations, the number of
constraints, the lower and upper bounds of the box constraints, the values of
the relaxation parameters, the initial values of the steering parameters and
the steering sequence. In all our experiments, the steering sequence of
Algorithm \ref{alg:steer} assumed the form
\begin{equation}
\sigma_{k}=\frac{\sigma}{k+1} \label{eq:sigma}%
\end{equation}
with a fixed user-chosen constant $\sigma.$ The main performance measure is
the value of $f(x^{k})$, plotted as a function of the iteration index $k$.

\subsection{Test problem description}

There are three types of constraints in our test problems: Box constraints,
linear constraints and quadratic constraints. Some of the numerical values
used to generate the constraints are uniformly distributed random numbers,
lying in the interval $\tau=[\tau_{1},\tau_{2}]$, where $\tau_{1}$ and
$\tau_{2}$ are user-chosen pre-determined values.

The $n$ box constraints are defined by%
\begin{equation}
\ell_{j}\leq x_{j}\leq u_{j},\quad j=1,2,\ldots,n
\end{equation}
where $\ell_{j},u_{j}\in\tau$ are the lower and upper bounds, respectively.
Each of the $N_{q}$ quadratic constrains is generated according to%
\begin{equation}
G_{i}(x)=\langle x,U_{i}x\rangle+\langle v^{i},x\rangle+\beta_{i},\quad
i=1,2,\ldots,N_{q}.
\end{equation}
Here $U_{i}$ is are $n\times n$ matrices defined by%
\begin{equation}
U_{i}=W_{i}\Lambda_{i}W_{i}^{T},
\end{equation}
the $n\times n$ matrices $\Lambda_{i}$ are diagonal, positive definite, given
by%
\begin{equation}
\Lambda_{i}=\operatorname*{diag}\left(  \delta_{1}^{i},\,\delta_{2}^{i}%
,\ldots,\delta_{n}^{i}\right)
\end{equation}
where $0<\delta_{1}^{i}\leq\,\delta_{2}^{i}\leq\ldots\leq\delta_{n}^{i}\in
\tau$ are generated randomly. The matrices $W_{i}$ are generated by
orthonormalizing an $n\times n$ random matrix, whose entries lie in the
interval $\tau$. Finally, the vector $v^{i}\in R^{n}$ is constructed so that
all its components lie in the interval $\tau$ and similarly the scalar
$\beta_{i}\in\tau$. The $N_{\ell}$ linear constraints are constructed in a
similar manner according to%
\begin{equation}
L_{i}(x)=\langle y^{i},x\rangle+\gamma_{i},\quad i=1,2,\ldots,N_{\ell}.
\end{equation}
Thus, the total number of constraints is $n+N_{q}+N_{\ell}$.

Table \ref{tests} summarizes the test cases used to evaluate and compare the
performance of Algorithms \ref{alg:ssp}, \ref{alg:steer} and \ref{alg:dan}. In
these eight experiments, we modified the value of the constant $\sigma$ in
(\ref{eq:sigma}), the interval $\tau$, the number of constraints, the number
of iterations, and the relative tolerance $\varepsilon$, used as a termination
criterion between subsequent iterations.

\begin{table}[h]
\begin{center}%
\begin{tabular}
[c]{|c||c|c|c|c|c|c|c|}\hline
Case & $\alpha$/$\sigma$/$\lambda$ & $\tau$ & $n$ & $N_{q}$ & $N_{\ell}$ &
Iterations & $\varepsilon$\\\hline\hline
1 & 1.1 &  & 3 & 5 & 5 & 1,000 & \\\hline
2 & 1.1 &  & 3 & 5 & 5 & 1,000 & \\\hline
3 & 1.98 &  & 3 & 5 & 5 & 1,000 & \\\hline
4 & 1.98 & $[-0.1,0.1]$ & 30 & 50 & 50 & 1,000 & 0.1\\\hline
5 & 1.98 & $[-10,10]$ & 30 & 50 & 50 & 100,000 & 0.1\\\hline
6 & 2 & $[-0.1,0.1]$ & 30 & 50 & 50 & 1,000 & 0.1\\\hline
7 & 3 & $[-10,10]$ & 3 & 5 & 5 & 1,000 & 0.1\\\hline
8 & 5 & $[-0.1,0.1]$ & 3 & 5 & 5 & 1,000 & 0.1\\\hline
\end{tabular}
\label{tests}
\end{center}
\caption{Test cases for performance evaluation}%
\end{table}In Table \ref{tests}, Cases 1 and 2 represent small-scale problems,
with a total of 13 constraints, whereas Cases 4--6 represent mid-scale
problems, with a total of 130 constraints. Cases 6--8 examine the case of over
relaxation, wherein the initial steering (relaxation) parameter is at least 2.

\subsection{Results}

\begin{center}
\begin{figure}[ptb]
\subfigure[Case 1]{\includegraphics[width=0.6\linewidth]{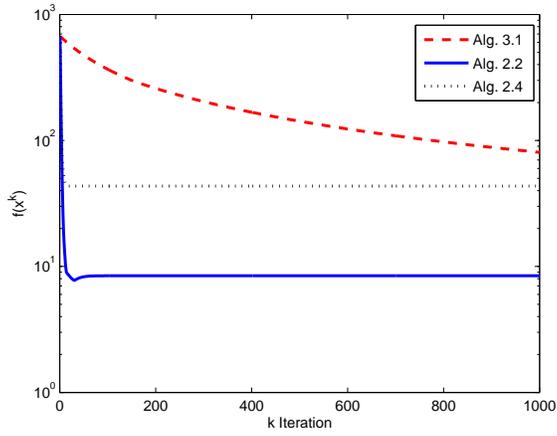}}
\subfigure[Case 2]{\includegraphics[width=0.6\linewidth]{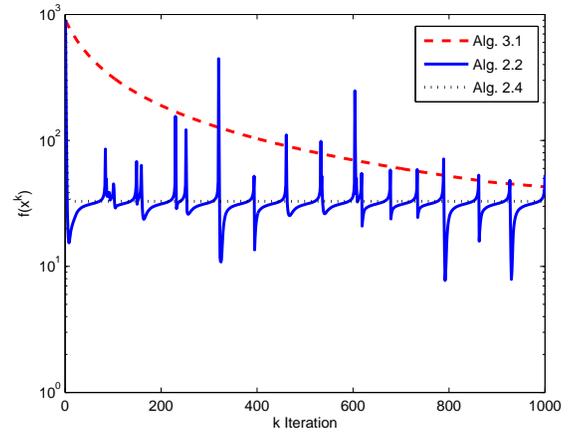}}
\subfigure[Case 3]{\includegraphics[width=0.6\linewidth]{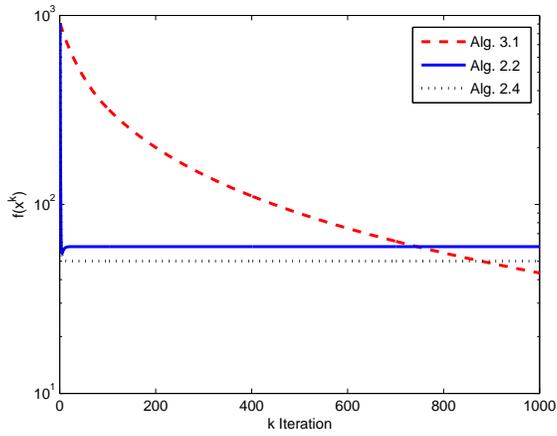}}\caption{Simulation
results for a small-scale problem, comparing Algorithms \ref{alg:ssp},
\ref{alg:steer} and \ref{alg:dan}. }%
\label{case1}%
\end{figure}
\end{center}
\begin{figure}[ptb]
\subfigure[Case 6]{\includegraphics[width=0.6\linewidth]{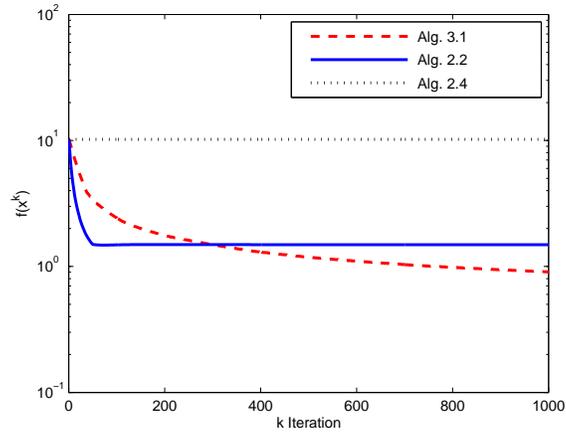}}
\subfigure[Case 7]{\includegraphics[width=0.6\linewidth]{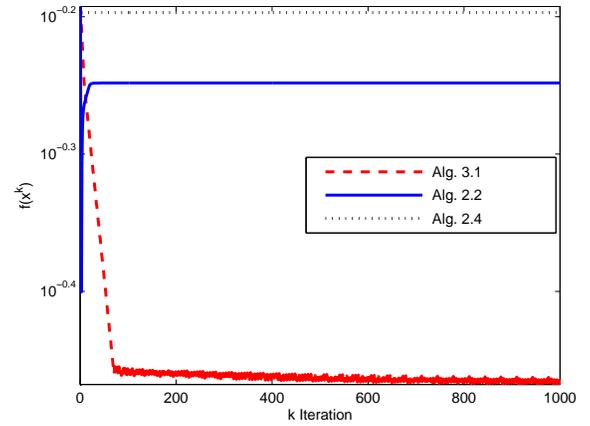}}
\subfigure[Case 8]{\includegraphics[width=0.6\linewidth]{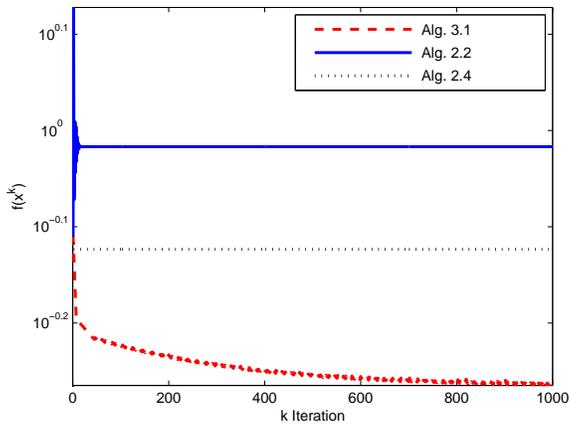}}\caption{Simulation
results for small- and mid-scale problems with overrelaxation, comparing
Algorithms \ref{alg:ssp}, \ref{alg:steer} and \ref{alg:dan}. }%
\label{case6}%
\end{figure}

The results of our experiments are depicted in Figures~\ref{case1}%
--\ref{case6}. The results of Cases 1--3 are shown in Figures~\ref{case1}%
(a)--\ref{case1}(c), respectively. It is seen that in Case 1 Algorithm
\ref{alg:ssp} has better initial convergence than Algorithms \ref{alg:steer}
and \ref{alg:dan}. However, in Case 2, Algorithm \ref{alg:steer} yields fast
and smooth initial behavior, while Algorithm \ref{alg:ssp} oscillates
chaotically. Algorithm \ref{alg:dan} exhibits slow initial convergence,
similarly to Case 1. In Case 3, Algorithm \ref{alg:dan} supersedes the
performance of the other two algorithm, since it continues to converge toward
zero. However, none of the algorithms detects a feasible solution, since none
converged to the tolerance threshold after the maximum number of iterations.

\begin{figure}[ptb]
\subfigure[Case 4]{\includegraphics[width=0.6\linewidth]{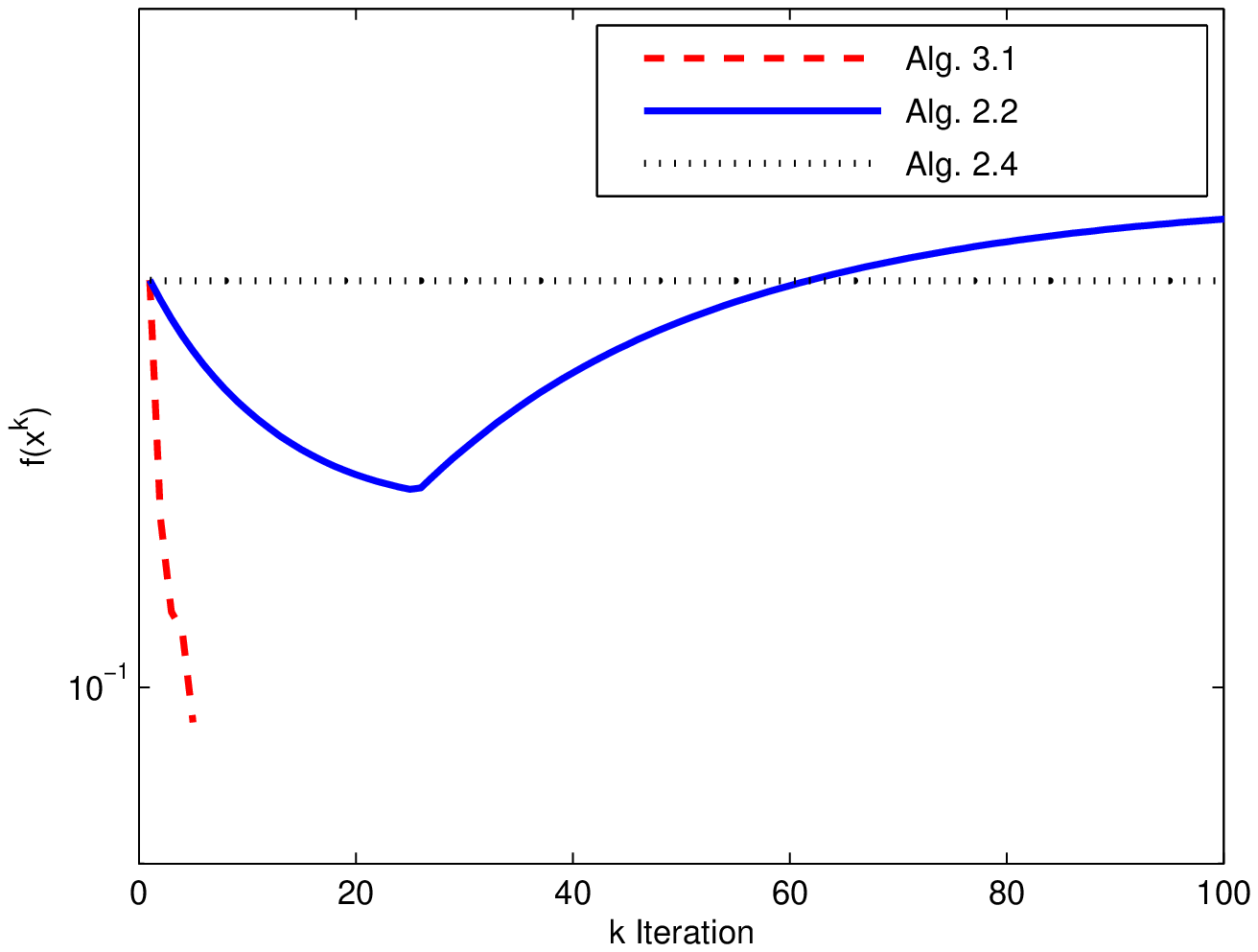}}
\subfigure[Case 5]{\includegraphics[width=0.6\linewidth]{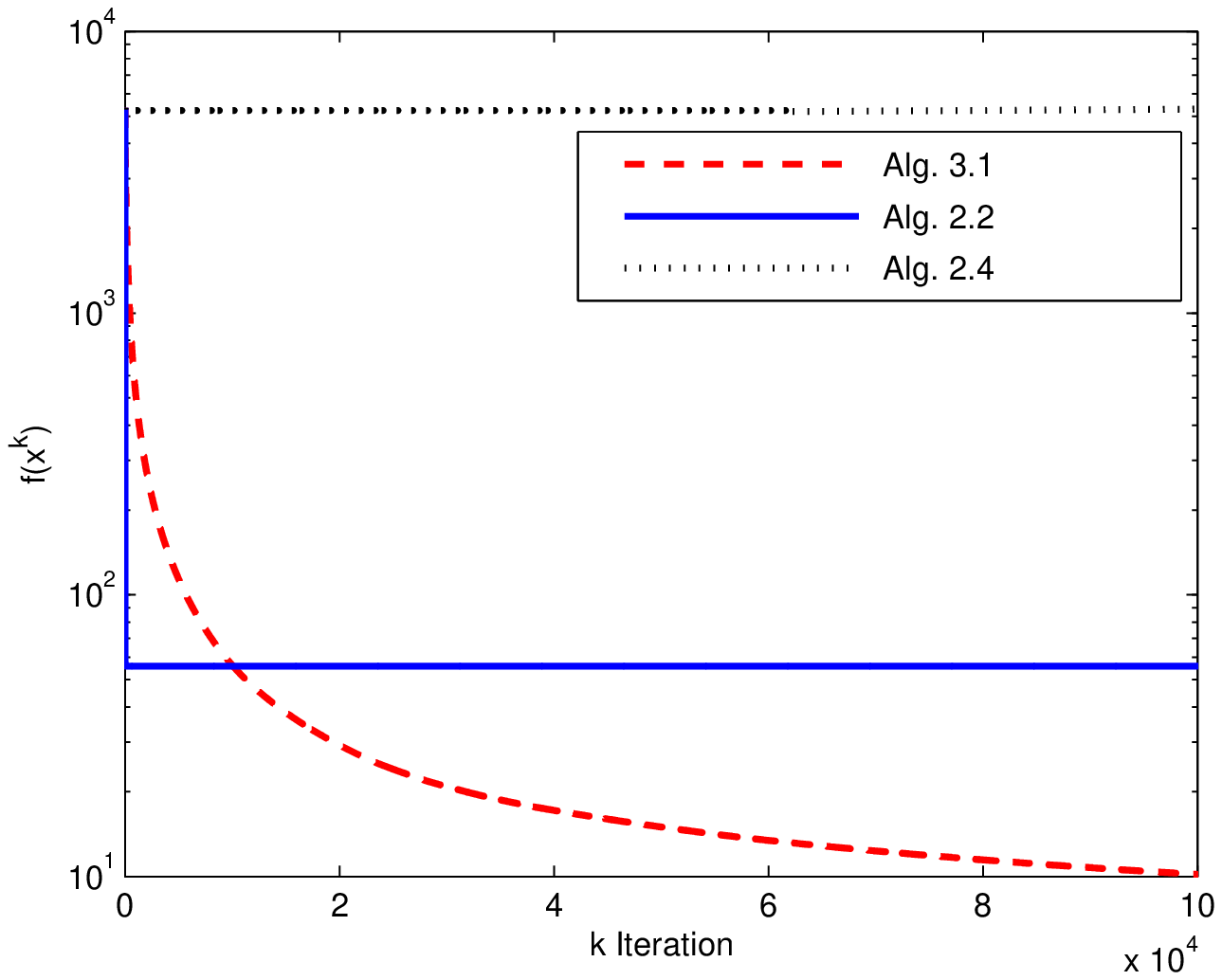}}\caption{Simulation
results for a mid-scale problem, comparing Algorithms \ref{alg:ssp},
\ref{alg:steer} and \ref{alg:dan}. }%
\label{case4}%
\end{figure}

The mid-sized problems of Cases 4 and 5 are depicted by Figures~\ref{case4}(a)
and \ref{case4}(b). Figure~\ref{case4}(a) shows that Algorithm \ref{alg:dan}
detects a feasible solution, while both Algorithms \ref{alg:ssp} and
\ref{alg:steer} fail to detect such a solution. The curve of Algorithm
\ref{alg:dan} in Fig. \ref{case4}(a) stops when it reaches the feasible point
detection tolerance, which is 0.1. Once the point is detected, there is no
need to further iterate, and the process stops. The curve for Algorithm
\ref{alg:ssp} in this figure shows irregular behavior since it searches for a
feasible solution without reaching the detection threshold of 0.1 and
accumulated numerical errors start to affect it. Figure~\ref{case4}(b) shows a
phenomenon similar to the one observed in the small-scale problem: Algorithm
\ref{alg:dan} continues to seek for a feasible solution, while Algorithms
\ref{alg:ssp} and \ref{alg:steer} converge to a steady-state, indicating
failure to detect a feasible solution.

In the experiments, Cases 6--8, Algorithm \ref{alg:dan} outperforms the other
algorithms, arriving very close to finding feasible solutions. It should be
observed that the behavior of Algorithm \ref{alg:dan} observed above is the
result of the way in which the relaxation parameters $\lambda_{k}$ are
self-regulating their sizes. In Algorithm \ref{alg:dan} the relaxation
parameter $\lambda_{k}$ can be chosen (see Equation (3.6)) to be any number of
the form%
\begin{equation}
\lambda_{k}=\beta_{k}\frac{\max(0,f(x^{k}))}{M^{2}}+2(1-\beta_{k})\frac
{\max(0,f(x^{k}))}{M^{2}}=(2-\beta_{k})\frac{\max(0,f(x^{k}))}{M^{2}},
\end{equation}
where $\beta_{k}$ runs over the interval $\left[  0,1\right]  .$ Consequently,
the size of $\lambda_{k}$ can be very close to zero when $x^{k}$ is close to a
feasible solution (no matter how $\beta_{k}$ is chosen in $\left[  0,1\right]
)$. Also, $\lambda_{k}$ may happen to be much larger then $2$ when $x^{k}$ is
far from a feasible solution and the number $f(x^{k})$ is large enough (note
that $2-\beta_{k}$ stays between 1 and 2). So, Algorithm \ref{alg:dan} is
naturally under- or over- relaxing the computational process according to the
relative position of the current iterate $x^{k}$ to the feasibility set of the
problem. As our experiments show, in some circumstances, this makes Algorithm
\ref{alg:dan} behave better then the other procedures we compare it with. At
the same time, the self-regulation of the relaxation parameters, which is
essential in Algorithm \ref{alg:dan}, may happen to reduce the initial speed
of convergence of this procedure, that is, Algorithm \ref{alg:dan} may require
more computational steps in order to reach a point $x^{k}$ which is close
enough to the feasibility set such that its self-regulatory features to be
really advantageous for providing a very precise solution of the given problem
(which the other procedures may fail to do since they may became stationary in
the vicinity of the feasibility set). Another interesting feature of Algorithm
\ref{alg:dan}, which differentiates it from the other algorithms we compare it
with, is its essentially non-simultaneous character: Algorithm \ref{alg:dan}
does not necessarily ask for $w_{i}^{k}>0$ for all $i\in\{1,...,m\}.$ The set
of positive weights $w_{i}^{k}$ which condition the progress of the algorithm
at step $k$ essentially depends on the current iterate $x^{k}$ (see (3.4)) and
allows reducing the number of subgradients needed to be computed at each
iterative step (in fact, one can content himself with only one $w_{i}^{k}>0$
and, thus, with a single subgradient $\xi_{i}^{k}$). This may be advantageous
in cases when computing subgradients is difficult and, therefore, time
consuming.\bigskip

The main observations can be summarized as follows:

\begin{enumerate}
\item Algorithm \ref{alg:dan} exhibits faster initial convergence than the
other algorithms in the vicinity of points with very small $f(x^{k})$. When
the algorithms reach points with small $f(x^{k})$ values, then Algorithm
\ref{alg:dan} tends to further reduce the value of $f(x^{k})$, while the other
algorithms tend to converge onto a constant steady-state value.

\item The problem dimensions in our experiments have little impact on the
behavior of the algorithms.

\item All the examined small-scale problems have no feasible solutions. This
can be seen from the fact that all three algorithms stabilize around
$f(x^{k})=50$.

\item The chaotic oscillations of Algorithm \ref{alg:ssp} in the underrelaxed
case is due to the fact that this algorithm has no internal mechanism to
self-adapt its progress to the distance between the current iterates and the
sets whose intersections are to be found. This phenomenon can hardly happen in
Algorithm \ref{alg:dan} because its relaxation parameters are self-adapting to
the size of the current difference between successive iterations. This is an
important feature of this algorithm. However, this feature also renders it
somewhat slower than the other algorithms.

\item In some cases, Algorithms \ref{alg:ssp} and \ref{alg:steer} indicate
that the problem has no solution. In contrast, Algorithm \ref{alg:dan}
continues to make progress and seems to indicate that the problem has a
feasible solution. This phenomenon is again due to the self-adaptation
mechanism, and can be interpreted in one of the following ways: (a) The
problem indeed has a solution but Algorithms \ref{alg:ssp} and \ref{alg:steer}
are unable to detect it (because they stabilize too fast). Algorithm
\ref{alg:dan} detects a solution provided that it is given enough running
time; (b) The problem has no solution and then Algorithm \ref{alg:dan} will
stabilize close to zero, indicating that the problem has no solution, but this
may be due to computing (round-off) errors. Thus, a very small perturbation of
the functions involved in the problem may render the problem feasible.
\end{enumerate}

\section{Conclusions}

We have studied here mathematically and experimentally subgradient projections
methods for the convex feasibility problem. The behavior of the fully
simultaneous subgradient projections method in the inconsistent case is not
known. Therefore, we studied and tested two options. One is the use of
steering parameters instead of relaxation parameters and the other is a
variable relaxation strategy which is self-adapting. Our small-scale and
mid-scale experiments are not decisive in all aspects and call for further
research. But one general feature of the algorithm with the self-adapting
strategical relaxation is its stability (non-oscillatory) behavior and its
relentless improvement of the iterations towards a solution in all cases. At
this time we have not yet refined enough our experimental setup. For example,
by the iteration index $k$ on the horizontal axes of our plots we consider a
whole sweep through all sets of the convex feasibility problem, regardless of
the algorithm. This is a good first approximation by which to compare the
different algorithms. More accurate comparisons should use actual run times.
Also, several numerical questions still remain unanswered in this report.
These include the effect of various values of the constant $\sigma$ as well as
algorithmic behavior for higher iteration indices. In light of the
applications mentioned in Section \ref{sect:intro}, higher dimensional
problems must be included. These and other computational questions are
currently investigated. \bigskip

\textbf{Acknowledgments. }We gratefully acknowledge the constructive comments
of two anonymous referees which helped us to improve an earlier version of
this paper. This work was supported by grant No.~2003275 of the United
States-Israel Binational Science Foundation (BSF), by a National Institutes of
Health (NIH) grant No. HL70472, by The Technion - University of Haifa Joint
Research Fund and by grant No.~522/04 of the Israel Science Foundation (ISF)
at the Center for Computational Mathematics and Scientific Computation (CCMSC)
in the University of Haifa.

\end{document}